\documentclass[12pt,a4paper]{article}
\usepackage{amssymb}

 \newcommand{\bm}[1]{\mbox{\boldmath $#1$}}

 \newcommand{\NN}{\mathbb N}
 \newcommand{\RR}{\mathbb R}
 
 \newcommand{\ZZ}{\mathbb Z}
 \newcommand{\II}{\mathbb I}
 \newcommand{\JJ}{\mathbb J}
 \newcommand{\HH}{\mathbb H}
 \newcommand{\QQ}{\mathbb Q}
 \newcommand{\KK}{\mathbb K}
 
 \newcommand{\kk}{{K}}
 \newcommand{\oo}{{\scriptstyle {\cal O}}}
 \newcommand{\OO}{{\cal O}}
 \newcommand{\LL}{{\cal L}}
 \newcommand{\MM}{{\cal M}}
 \newcommand{\aaa}{\mathfrak a}
 \newcommand{\ppp}{\mathfrak p}
 \newcommand{\PPP}{\mathfrak P}
 
 \newcommand{\rf}{\mbox{\small $\sqrt{5}$}}
 \newcommand{\rz}{\mbox{\small $\sqrt{2}$}}

\newtheorem{defin}{Definition}
\newtheorem{theorem}{Theorem}
\newtheorem{lemma}{Lemma}
\newtheorem{prop}{Proposition}
\newtheorem{coro}{Corollary}

 \newcommand{\be}{\begin{equation}}
 \newcommand{\ee}{\end{equation}}

\begin{document}
 \bibliographystyle{unsrt}

 \parindent0pt

\begin{center}
{\large \bf Similarity submodules and root systems in four dimensions}

 \end{center}
 \vspace{3mm}
 
 \begin{center}  {\sc Michael Baake}$^{\, \rm a,}$\footnote{Heisenberg-Fellow} 
             and {\sc Robert V. Moody}$^{\, \rm b}$ 
 
\vspace{5mm}

 a) Institut f\"ur Theoretische Physik, Universit\"at T\"ubingen, \\

 Auf der Morgenstelle 14, D-72076 T\"ubingen, Germany  \\

 b) Department of Mathematical Sciences, University of Alberta, \\

 Edmonton, Alberta T6G 2G1, Canada
 
\end{center}
 
\vspace{8mm}
\centerline{Dedicated to H.~S.~M.~Coxeter}
\vspace{5mm}
 
\begin{abstract} 
Lattices and $\ZZ$-modules in Euclidean space possess an infinitude of
subsets that are images of the original set under similarity transformation.
We classify such self-similar images according to their indices for
certain 4D examples that are related to 4D root systems, both
crystallographic and non-crystallographic. We encapsulate their statistics 
in terms of Dirichlet series generating functions
and derive some of their asymptotic properties.
\end{abstract}

\parindent15pt

\subsection*{Introduction}

This paper begins with the problem of determining the self-similar images
of certain lattices and $\ZZ$-modules in four dimensions and ends in the
enchanting garden of Coxeter groups, the arithmetic of several quaternion
rings, and the asymptotics of their associated zeta functions. The main
results appear in Theorems \ref{thm-lattices} and \ref{thm-modules}.

The symmetries of crystals are of fundamental physical importance and,
along with the symmetries of lattices, have been studied by
mathematicians, crystallographers and physicists for ages. 
The recent interest in
{\em quasicrystals}, which are non-crystallographic yet still highly
ordered structures, has naturally led to speculation about the role of
symmetry in this new context. Here, however, it is apparent that a
different set of symmetry concepts is appropriate, notably because
translational symmetry is either entirely lacking or at least
considerably restricted in scope.

One of the most obvious features of quasicrystals is their tendency to
have copious inflationary self-similarity.  Thus instead of groups of
isometries, one is led to semi-groups of self-similarities that map
a given (infinite) point set into an inflated copy lying within
itself. Ordinary point symmetries then show up as a small part
of this, namely as the ``units'', i.e. as the (maximal
subgroups of) invertible elements.

The importance of self-similarities is well-known, and has also been used to gain 
insight into the colour symmetries of crystals \cite{Schwarz1,Schwarz2} and, more 
recently, of quasicrystals \cite{Baake97,BM1,Ron}. The focus in the latter cases 
was on highly symmetric examples in the plane and in 3-space since they are 
obviously of greatest physical importance. 

In this contribution, we extend the investigation of self-similarities to certain 
exceptional examples in 4-space, namely the hypercubic lattices (of which there are 
two, represented by the primitive hypercubic lattice $\ZZ^4$ and the root lattice 
$D^{}_4$ or, equivalently, its weight lattice $D^*_4$) and the icosian ring, seen as 
the $\ZZ$-span of the root system of the non-crystallographic Coxeter group $H_4$, 
see \cite{CS,Humphreys,Cox80} for notation and background material. The case of
lattices has recently also been investigated by Conway, Rains and Sloane \cite{CRS}
who are generally interested in the question under which conditions similarity
sublattices of a given index exist. Their methods are complementary to ours,
and more general, but do not seem to give direct access to the full combinatorial
problem which we can solve here.

It is useful to digress briefly to discuss the role of $H_4$ in the context of 
aperiodic order. It is a remarkable fact that the non-crystallographic point 
symmetries relevant to essentially all known physical quasicrystals are actually 
Coxeter groups, namely the dihedral groups $I_2(k)$, $k = 5,8,10,12$ and the 
icosahedral group $H_3$ \cite{BJKS}. Apart from the remaining infinite series of
dihedral groups, the only indecomposable non-crystallographic Coxeter group is 
$H_4$ of order 14400, which is the most interesting of them all. In spite of being 
four-dimensional in nature, there are several good physical reasons to explore the 
symmetry expressed by this group, and because of its connections with exceptional 
objects in mathematics, including the root lattice $E_8$, there are good 
mathematical reasons, too.

{}From the point of view of quasicrystals, $H_4$ may be viewed as the top member of 
the series $H_2 := I_2(5) \subset H_3 \subset H_4$ of which the first two have been 
the subject of great attention, see \cite{Kramer} and references therein, while $H_4$ 
appears as symmetry group of the Elser-Sloane quasicrystal \cite{ES}. Mathematically, 
this family belongs together and $H_4$ is the natural parent of the others. 

Now, $H_4$ has a natural quaternionic interpretation which arises
as follows. The group of norm $1$ units of the real quaternion algebra
is easily identifiable with $\mbox{SU}(2)$, see \cite{Koecher} for
background material and notation. Using the $2$-fold cover of 
$\mbox{SO}(3)$ by $\mbox{SU}(2)$ we can find (in many ways) a 
$2$-fold cover of the icosahedral group inside $\mbox{SU}(2)$. 
This is the binary icosahedral group $I$ of order $120$ \cite[p.\ 69]{CM}. 
The point set $I$ is a beautiful object, namely the set of
vertices of the exceptional regular $4$-polytope called the 
$600$-cell\footnote{The 600-cell and its dual, the 120-cell, are
two of the three exceptional regular polytopes in 4-space \cite[p.\ 292]{Cox73}. 
The remaining one, the 24-cell, also occurs, later in this paper.
Beyond 4 dimensions, the only regular polytopes are the simplices,
the hypercubes, and their duals, the cross-polytopes (or
hyperoctahedra).} \cite[Ch.\ 22]{Cox80} (under the standard topology of 
$\RR^4$ carried by the quaternion algebra).
The $\ZZ$-span of $I$ is a ring $\II$, dubbed by Conway the
ring of {\em icosians}. This ring is closed under complex conjugation 
and under left and right multiplications by elements of $I$. The group
of symmetries of $\II$ so obtained is isomorphic to $H_4$
acting as a reflection group in $\RR^4$.  The ring
$\II$ is itself quite a remarkable object. It is naturally 
a rank $4$-module over $\ZZ[\tau]$, $\tau := (1 + \sqrt 5)/2$, and 
a rank $8$-module over $\ZZ$ (so it is certainly dense in the ambient
space $\RR^4$). In fact, as an aside, it has a canonical interpretation 
as the root lattice of type  $E_8$ with $I \cup \tau I$ making
up the $240$ roots of $E_8$. 
Restricting to the pure quaternions
puts us in the $3$-dimensional icosahedral case, and by
further restriction we can get the $H_2$ situation.

Now we can state our problem for the icosian case. We have 
pointed out that the additve group $\II$ has a large finite
group of rotational symmetries coming from the left and right
multiplications by elements of $I$ -- namely $120^2/2=7200$ such 
symmetries. But we have seen that in the study of quasicrystals
we have to pay attention to self-similarities, too. So we are now 
also interested in rotation-inflations of
$\RR^4$ that map $\II$ into, but not necessarily onto, itself. 
Each such self-similarity maps $\II$ onto some submodule of finite
index, and our question is to determine these images and to count
the number of different similarity submodules of a given index.
This leads us to introduce a suitable Dirichlet series
generating function which encodes the counting information 
and its asymptotic properties all at once, and indeed
determining its exact form is a number theoretical problem
that depends heavily on the fact that $\II$ can be interpreted
as a maximal order in the split quaternionic algebra over the
quadratic field $\QQ(\rf)$. 

The other situation that we wish to discuss in this paper
is crystallographic in origin, but it nevertheless is very much
the same problem. It is well known that there are two hypercubic
lattices in $4$-space \cite{Schwarz1,Brown}, namely the primitive and
the centred one (face-centred and body-centred are equivalent in $4$-space
by a similarity tranformation). Let us take $\ZZ^4$ and the root lattice
$D_4$ as suitable representatives. Note that they have different holohedries,
namely one of order 1152 (denoted by 33/16 in \cite[Fig.\ 7]{Brown}) for
$D_4$, which coincides with the automorphism group of the underlying root system,
and an index 3 subgroup (denoted by 32/21 in \cite[Fig.\ 7]{Brown}) for $\ZZ^4$.
Due to the previous remark, we may
take the weight lattice $D^*_4$ instead of $D_4$ if we wish, and we will
frequently do so. Given any of these cases, we want to determine
the Dirichlet series generating function for the sublattices that are 
self-similar images of it. What makes this situation tractable is that,
parallel to the icosian case, there is a highly structured
algebraic and arithmetic object in the background, namely Hurwitz'
ring of integral quaternions \cite{Hurwitz,Cox46}. In our setting,
it is $\JJ=D^*_4$, and it is again a maximal order, this time of
the quaternionic algebra over the rationals, $\QQ$. The results
for the Hurwitzian and icosian cases are striking in their
similarity.
Another example is that of the maximal order in $\HH(\QQ(\rz))$ the
treatment of which equals that of $\II$ whence we only state the results.

The structure of the paper is as follows. In the next Section, we
set the scene by collecting some methods and results from algebra,
and algebraic number theory in particular. This will be done in
slightly greater detail than necessary for a mathematical audience,
but since there is also considerable interest in this type of problem
from the physics community, we wish to make the article more
self-contained and readable this way. The two following Sections give
the results, first for lattices, and then for modules.
We close with a brief discussion of related aspects and provide
an Appendix with material on the asymptotics of arithmetic functions
defined through Dirichlet series.

\subsection*{Preliminaries and Recollections}

We shall need a number of results from algebraic number theory,
both commutative and non-commutative. {}First of all, we need, of course,
the arithmetic of $\ZZ$, the ring of integers in the field $\QQ$.
All ideals of $\ZZ$ are principal, and they are of the form
$\aaa = m\ZZ$ with $m\in\ZZ$. If $\aaa\neq 0$, the index is
$[\ZZ:\aaa]=|m|$. The corresponding zeta function, which can be seen as 
the Dirichlet series generating function for the number of non-zero ideals 
of a given index, is Riemann's zeta function itself \cite{Apostol}
\be \label{Riemann}
    \zeta(s) \; = \; \sum_{\aaa\subset\ZZ}\frac{1}{[\ZZ:\aaa]^s}
             \; = \; \sum_{m=1}^{\infty} \frac{1}{m^s}
             \; = \; \prod_{p\in{\cal P}}
                     \frac{1}{1-p^{-s}} \, .
\ee
Here, ${\cal P}$ denotes the set of (rational) primes, and the
second representation of the Riemann zeta function is its
{\em Euler product expansion}. It is possible because the number of
ideals of index $m$ is a multiplicative arithmetic function,
a situation that we shall encounter throughout the article.

Next, we need the analogous objects for the real quadratic field
$\QQ(\rf)=\QQ(\tau)$. The ring of integers turns out to be
\be
    \ZZ[\tau] \; = \; \{ m + n \tau \mid m,n \in \ZZ \}
\ee
where $\tau = (1+\rf)/2$ is the fundamental unit of $\ZZ[\tau]$,
i.e.\ all units are obtained as $\pm\tau^m$ with $m\in\ZZ$.
Again, $\ZZ[\tau]$ is a principal ideal domain, and hence a
unique factorization domain \cite[ch.\ 15.4]{Hardy}.
The zeta function is the Dedekind zeta function 
\cite[\S 11]{Zagier} defined by
\be
  \zeta_{\QQ(\tau)}^{}(s) \; = \; 
    \sum_{\aaa\subset\ZZ[\tau]} \frac{1}{[\ZZ[\tau]:\aaa]^s}
    \; = \; \sum_{m=1}^{\infty} \frac{a(m)}{m^s}
\ee
where $\aaa$ runs through the non-zero ideals of $\ZZ[\tau]$ and
$[\ZZ[\tau]:\aaa]$ is the norm of $\aaa$.
If $\aaa = \alpha\ZZ[\tau]$, it is given by 
\be
     [\ZZ[\tau]:\aaa] \; = \; |\mbox{N}(\alpha)| \; = \;
     |\alpha \alpha'|
\ee
where $'$ denotes algebraic conjugation in $\QQ(\tau)$, defined by
$\tau \mapsto 1-\tau$.

Explicitly, the zeta function reads (see the Appendix for details):
\begin{eqnarray} \label{zeta2}
  \zeta_{\QQ(\tau)}^{}(s) &\!\!\! = &\!\!\!
            \frac{1}{1-5^{-s}} \cdot
            \prod_{p \equiv \pm 1 \; (5)} 
                   \frac{1}{(1-p^{-s})^2} \cdot
            \prod_{p \equiv \pm 2 \; (5)} 
                   \frac{1}{1-p^{-2s}}  \\
&\!\!\! = &\!\!\! \mbox{\small $
     1+\frac{1}{4^s}+\frac{1}{5^s}+\frac{1}{9^s}+\frac{2}{11^s}+
       \frac{1}{16^s}+\frac{2}{19^s}+\frac{1}{20^s}+\frac{1}{25^s}+
       \frac{2}{29^s}+\frac{2}{31^s}+
       \frac{1}{36^s}+\frac{2}{41^s}  
     + \cdots $}  \nonumber
\end{eqnarray}

As before, $a(m)$ is a multiplicative arithmetic function, i.e.\
$a(mn)=a(m)a(n)$ for coprime $m,n$.
It is thus completely specified by its value for $m$ being a
prime power, and from the Euler product in (\ref{zeta2}) one
quickly derives that $a(5^r)=1$ (for $r\geq 0$). Then, for
primes $p\equiv\pm 2$ $(5)$, one obtains $a(p^{2r+1})=0$ and
$a(p^{2r})=1$, while for primes $p\equiv\pm 1$ $(5)$, the
result is $a(p^r) = r+1$. 

One benefit of relating the numbers $a(m)$ to zeta functions 
with well-defined analytic behaviour is that one can rather easily
determine the asymptotic behaviour of $a(m)$ from the poles of 
the zeta function, see the Appendix for a summary. In this case, the
function $a(m)$ is constant on average, the constant being the
residue of $\zeta^{}_{\QQ(\tau)}(s)$ at its right-most pole, $s=1$.
Explicitly, we get
\be \label{asymp1}
    \lim_{N\to\infty}
    \frac{1}{N} \sum_{m=1}^{N} a(m) \; = \;
    \mbox{res}_{s=1} \, \zeta^{}_{\QQ(\tau)}(s)
    \; = \; \frac{2\log(\tau)}{\sqrt{5}} 
    \; \simeq \; 0.430409 \, .
\ee

We shall also need the zeta function of the quadratic field $\QQ(\rz)$, where
$\ZZ[\rz]$ is the corresponding ring of integers, and $1+\rz$ its fundamantal unit
\cite{Hardy}. The zeta function reads
\begin{eqnarray} \label{zeta3}
  \zeta_{\QQ({\scriptscriptstyle \sqrt{2}})}^{}(s) &\!\!\! = &\!\!\!
            \frac{1}{1-2^{-s}} \cdot
            \prod_{p \equiv \pm 1 \; (8)} 
                   \frac{1}{(1-p^{-s})^2} \cdot
            \prod_{p \equiv \pm 3 \; (8)} 
                   \frac{1}{1-p^{-2s}}  \\
&\!\!\! = &\!\!\! \mbox{\small $
     1+\frac{1}{2^s}+\frac{1}{4^s}+\frac{2}{7^s}+\frac{1}{8^s}+
       \frac{1}{9^s}+\frac{2}{14^s}+\frac{1}{16^s}+\frac{2}{17^s}+
       \frac{1}{18^s}+\frac{2}{23^s}+
       \frac{1}{25^s}+\frac{2}{28^s}  
     + \cdots $}  \nonumber
\end{eqnarray}
The coefficients are $a(2^r)=1$, $a(p^r)=r+1$ for $p\equiv \pm 1$ (8)
and $a(p^r)=0$ resp.\ $1$ for $p\equiv \pm 3$ (8) and $r$ odd resp.\ even.
The asymptotic behaviour is given by
$\lim^{}_{N\to\infty} {1\over N}\sum_{m\leq N} a(m) =
\log(1+\rz)/\rz \simeq 0.623225$.

Let us now move to the non-commutative results we shall need. We will 
be concerned with the quaternionic algebra $\HH(\kk)$, mainly 
over the field $\kk=\QQ$ or over $\kk=\QQ(\tau)$. The case of $\kk=\QQ(\rz)$
is treated more as an aside. In all cases, we are
interested in the  corresponding ring of integers, these being the 
{\em Hurwitzian ring}  $\,\JJ$, the {\em icosian ring} $\,\II$,
and the {\em cubian ring} $\,\KK$.
These are maximal orders in their respective quaternionic
algebras \cite{Hurwitz,Vigneras}.

Let us first consider $\JJ=D^*_4$. In terms of the standard 
basis\footnote{The defining relations \cite{Cox46} are:
$\bm{i}^2=\bm{j}^2=\bm{k}^2=\bm{ijk}=-1$.}
$1,\bm{i},\bm{j},\bm{k}$ of $\HH(\QQ)$, $\JJ$ consists of the points
$(x_0,x_1,x_2,x_3)$ whose coordinates $x_i$ either all lie in 
$\ZZ$ or all in $\ZZ+\frac{1}{2}$. Though non-commutative, $\JJ$ is
still a principal ideal domain, i.e.\ all left-ideals (and also all
right-ideals) are principal \cite{Hurwitz}. Consequently, we 
have unique factorization up to units, the units being the 24 elements
obtained from $(\pm 1,0,0,0)$ plus permutations and 
from ${1\over 2}(\pm 1,\pm 1,\pm 1,\pm 1)$. 
These 24 units form a group, $\JJ^{\times}$, that is isomorphic to the binary 
tetrahedral group \cite[p.\ 69]{CM}. The number of non-zero 
left-ideals (or, equivalently, of non-zero right-ideals) can 
now be counted explicitly, which results in the corresponding
zeta function $ \zeta^{}_{\JJ}(s) \; = \; \sum_{\aaa} 
\frac{1}{[\JJ:\aaa]^s}$, where $\aaa$ runs over the non-zero left ideals of 
$\JJ$, see \cite[Sec.\ VII, {\S} 8 and {\S} 9]{Deuring} for details on zeta 
functions of quaternionic algebras. The result is 
\cite[{\S} 63, A.~15]{Scheja}:
\begin{itemize}
\item  The zeta function of Hurwitz' ring of integer quaternions, $\JJ$,
       is given by\footnote{Note that the formula given in the first line
       after \cite[Eq.\ 3.17]{Baake} contains a misprint in the prefactor.}
\be \label{Hurwitz-zeta}
    \zeta^{}_{\JJ}(s) \; = \; (1-2^{1-2s}) \cdot 
                 \zeta(2s) \, \zeta(2s-1) \, .
\ee
\end{itemize}
Using (\ref{Riemann}), one can easily determine the first few terms
\be
  \zeta^{}_{\JJ}(s) \; = \;
  \mbox{\small $
     1+\frac{1}{4^s}+\frac{4}{9^s}+\frac{1}{16^s}+\frac{6}{25^s}+
       \frac{4}{36^s}+\frac{8}{49^s}+\frac{1}{64^s}+\frac{13}{81^s}+
       \frac{6}{100^s}+\frac{12}{121^s}+\frac{4}{144^s}+ 
       \cdots $ }
\ee
The possible indices of ideals are the squares of integers. It is
thus convenient to write the Dirichlet series as
\be \label{coeff1}
  \zeta^{}_{\JJ}(s) \; = \; \sum_{m=1}^{\infty}
          \frac{a^{}_{\JJ}(m)}{m^{2s}}
\ee
so that $a^{}_{\JJ}(m)$ is actually the number of left-ideals of
index $m^2$ (rather than $m$).
This then results in $a^{}_{\JJ}(2^r)=1$ (for $r\geq 0$) and in
$a^{}_{\JJ}(p^r)=(p^{r+1}-1)/(p-1)$ for odd primes. Let us add that
$a^{}_{\JJ}(m)$ is also the sum of the odd divisors of $m$, see
\cite[sequence M 3197]{SP} or \cite[sequence A 000593]{Sloane}.

Let us again briefly comment on the asymptotic behaviour. In this 
case, the average of $a^{}_{\JJ}(m)$ grows linearly with $m$, i.e.\
\be \label{asymp2}
     \frac{1}{N} \sum_{m=1}^{N} a^{}_{\JJ}(m) \; \sim \; 
     \frac{\pi^2}{24} N \quad \quad \quad (N \rightarrow\infty)
\ee
where the coefficient is half the residue of $\zeta^{}_{\JJ}(s)$
at its right-most pole, $s=1$. This can easily be calculated from
the details provided in the Appendix.

Note that the linear growth of the average of $a^{}_{\JJ}(m)$
stems from the definition used in (\ref{coeff1}). With the usual
definition (i.e.\ with denominators $m^s$ rather than $m^{2s}$),
the average would tend to a constant, as in (\ref{asymp1}).

Now, let us consider the analogous situation, with $\II$ being
a maximal order in the algebra $\HH(\QQ(\tau))$. Again, all
left-ideals (and all right-ideals) of $\II$ are principal, and
we also get unique factorization up to units again \cite{Vigneras}. 
The unit group $\II^{\times}$ of $\II$ consists of the 120 elements of 
the binary icosahedral group 
$I$ inside $\II$. Taking the unit quaternions $(1,0,0,0)$,
${1\over 2}(1,1,1,1)$, ${1\over 2}(\tau,1,-1/\tau,0)$ together with
all even permutations and arbitrary sign flips results in an explicit
choice of the group $I$, and hence of $\II$. Again defining the zeta 
function for one-sided ideals of $\II$ we have
\begin{itemize}
\item The zeta function of the icosian ring, $\II$, is
\be \label{icosian-zeta}
    \zeta^{}_{\II}(s) \; = \; 
                 \zeta^{}_{\QQ(\tau)}(2s) \, 
                 \zeta^{}_{\QQ(\tau)}(2s-1) \, .
\ee
\end{itemize}
This result follows from \cite[Ch.\ III, Prop.\ 2.1]{Vigneras}.
The first few terms of this series read
\be \label{ico-zeta-2}
  \zeta^{}_{\II}(s) \; = \;
  \mbox{\small $
     1+\frac{5}{16^s}+\frac{6}{25^s}+\frac{10}{81^s}+
       \frac{24}{121^s}+\frac{21}{256^s}+\frac{40}{361^s}+
       \frac{30}{400^s}+\frac{31}{625^s}+\frac{60}{841^s}+
       \frac{64}{961^s}+  \cdots $ }
\ee
The possible indices ($=$ denominators) are the squares of integers
that are representable by the quadratic form $x^2 + x y - y^2$,
i.e.\ of integers all of whose prime factors congruent to 2 or 3
(mod 5) occur with even exponent only. Using a definition analogous
to (\ref{coeff1}) above, the coefficient $a^{}_{\II}(m)$ is again a
multiplicative arithmetic function. It is given by
$a^{}_{\II}(5^r) = (5^{r+1}-1)/4$ (for $r\geq 0$), and, for primes
$p\equiv\pm 2$ $(5)$, by $a^{}_{\II}(p^{2r+1})=0$ and
$a^{}_{\II}(p^{2r})=(p^{2r+2}-1)/(p^2-1)$. {}Finally, for
$p\equiv\pm 1$ $(5)$, one finds $a^{}_{\II}(p^r) =
\sum_{l=0}^{r} (l+1) (r-l+1) p^l$. 
It is now listed as \cite[sequence A 035282]{Sloane}.

The asymptotic behaviour of $a^{}_{\II}(m)$ is similar to that of
$a^{}_{\JJ}(m)$ above, and we obtain
\be \label{asymp3}
   \frac{1}{N} \sum_{m=1}^{N} a^{}_{\II}(m) \; \sim \; 
     \frac{2\pi^4\log(\tau)}{375} N 
     \; \simeq \; 0.249997 \cdot N
     \quad \quad \quad (N \rightarrow\infty)
\ee
where the slope is again half the residue of $\zeta^{}_{\II}(s)$ at 
its right-most pole, $s=1$, see the Appendix for details.

Very similar is the situation of the ring $\KK$ in $\HH(\QQ(\rz))$,
generated as the $\ZZ[\rz]$-span of the basis
$\{1, (1+\bm{i})/\rz, (1+\bm{j})/\rz, (1+\bm{i}+\bm{j}+\bm{k})/2\}$.
The unit group $\KK^{\times}$ is the binary octahedral group of order 48
\cite[p.\ 69]{CM}, and
the symmetry group of $\KK$ contains that of $\JJ$ as an index 2 subgroup.
$\KK$ is again a maximal order and a principal ideal domain \cite{Vigneras}. 
The zeta function of $\KK$ reads \cite[Ch.\ III, Prop.\ 2.1]{Vigneras}
\begin{eqnarray} \label{cubian-zeta}
    \zeta^{}_{\KK}(s)  & \!\!\! =  &\!\!\!  
    \zeta^{}_{\QQ({\scriptscriptstyle \sqrt{2}})}(2s) \, 
    \zeta^{}_{\QQ({\scriptscriptstyle \sqrt{2}})}(2s-1) \\
&\!\!\! = &\!\!\! \mbox{\small $
     1+\frac{3}{4^s}+\frac{7}{16^s}+\frac{16}{49^s}+\frac{15}{64^s}+
       \frac{10}{81^s}+\frac{48}{196^s}+\frac{31}{256^s}+
       \frac{36}{289^s}+\frac{30}{324^s}+\frac{48}{529^s}+
       \frac{26}{625^s} + \cdots $}  \nonumber
\end{eqnarray}
and further details can be worked out in complete analogy to the icosian case.

Let us now briefly describe how the quaternions enter our (mainly
geometric) picture, and how they provide a parametrization of
$\mbox{(S)O}(4) = \mbox{(S)O}(4,\RR)$,
see \cite{Koecher} for details. The key is that
pairs of quaternions in $\HH(\RR)$, i.e.\ quaternions
$\bm{q}=(q^{}_0,q^{}_1,q^{}_2,q^{}_3)$ as written in the standard
basis $1,\bm{i},\bm{j},\bm{k}$ of the quaternion algebra over $\RR$,
induce an action on vectors of $\RR^4$ via
\be \label{action}
     M(\bm{q}^{}_1,\bm{q}^{}_2) \bm{x}^t \; = \;
       \bm{q}^{}_1 \bm{x} \, \overline{\bm{q}}^{}_2
\ee
where $M(\bm{q}^{}_1,\bm{q}^{}_2)\in \mbox{Mat}(4,\RR)$ and 
$\bm{x}^t$ is $\bm{x}$ written as a column vector in $\RR^4$. Evidently, 
for nonzero quaternions $\bm{q}^{}_1,\bm{q}^{}_2,\bm{r}^{}_1,\bm{r}^{}_2$ we have 
$M(\bm{q}^{}_1,\bm{q}^{}_2) = M(\bm{r}^{}_1,\bm{r}^{}_2)$ if and only if 
$\bm{r}^{}_1 = a \bm{q}^{}_1$ and $\bm{r}^{}_2 = a^{-1}\bm{q}^{}_2$ 
for some $a \in \RR \backslash \{0\}$.
With $\bm{q}^{}_1 = (a,b,c,d)$ and $\bm{q}^{}_2 = (t,u,v,w)$,
the matrix $M=M(\bm{q}^{}_1,\bm{q}^{}_2)$ reads explicitly
\be \label{matrix}
\left( \begin{array}{rrrr} 
     at\!+\!bu\!+\!cv\!+\!dw & \!-\!bt\!+\!au\!+\!dv\!-\!cw &
\!-\!ct\!-\!du\!+\!av\!+\!bw & \!-\!dt\!+\!cu\!-\!bv\!+\!aw \\
     bt\!-\!au\!+\!dv\!-\!cw &      at\!+\!bu\!-\!cv\!-\!dw &
\!-\!dt\!+\!cu\!+\!bv\!-\!aw &      ct\!+\!du\!+\!av\!+\!bw \\
     ct\!-\!du\!-\!av\!+\!bw &      dt\!+\!cu\!+\!bv\!+\!aw &
     at\!-\!bu\!+\!cv\!-\!dw & \!-\!bt\!-\!au\!+\!dv\!+\!cw \\
     dt\!+\!cu\!-\!bv\!-\!aw & \!-\!ct\!+\!du\!-\!av\!+\!bw &
     bt\!+\!au\!+\!dv\!+\!cw &      at\!-\!bu\!-\!cv\!+\!dw 
\end{array} \right)_.
\ee
It has determinant
\begin{eqnarray}
     \det(M(\bm{q}^{}_1,\bm{q}^{}_2)) & = &
      (a^2+b^2+c^2+d^2)^2\cdot(t^2+u^2+v^2+w^2)^2  \nonumber \\
     & = & \{(\bm{q}^{}_1)^2 \cdot (\bm{q}^{}_2)^2\}^2
\end{eqnarray}
and also fulfils
\be
     M M^t \; = \; \sqrt{\det{M}} \cdot 
                   \mbox{\large \bf 1}^{}_4 \, .
\ee

Consequently, when $|\bm{q}^{}_1| = |\bm{q}^{}_2| = 1$, we obtain a 4D
rotation matrix and the homomorphism 
$M: S^3\times S^3 \longrightarrow\mbox{SO}(4)$ 
provides the standard double cover of the rotation
group $\mbox{SO}(4)$ \cite{Koecher}, with 
$M(\bm{q}^{}_1,\bm{q}^{}_2) = M(-\bm{q}^{}_1,-\bm{q}^{}_2)$.
The orientation reversing transformations, i.e.\ the elements of 
$\mbox{O}(4)\backslash \mbox{SO}(4)$, are obtained by the mapping
$\bm{x} \mapsto \bm{q}^{}_1\,\overline{\bm{x}}\,\overline{\bm{q}}^{}_2$
with unit quaternions $\bm{q}^{}_1,\bm{q}^{}_2$.
Let us finally note that, for non-zero quaternions,
\be
   R(\bm{q}^{}_1,\bm{q}^{}_2) 
       \; := \; M(\frac{\bm{q}^{}_1}{|\bm{q}^{}_1|},
                  \frac{\bm{q}^{}_2}{|\bm{q}^{}_2|}) 
       \; = \; \frac{1}{|\bm{q}^{}_1 \bm{q}^{}_2|} \,
                M(\bm{q}^{}_1,\bm{q}^{}_2)
\ee
always gives a rotation matrix, which is handy for finding
suitable parametrizations of groups such as $\mbox{SO}(4,\QQ)$
or $\mbox{SO}(4,\QQ(\tau))$ in closely related problems, see 
\cite{Baake} for details.

\subsection*{Arguments in common}

In this Section, we focus on $\II$ versus $\JJ$ and carry the arguments 
as far as possible without having to separate the two rings $\II$ and $\JJ$ 
too seriously\footnote{The case $\OO=\KK$ is entirely parallel to that
of $\OO=\II$ and need not be spelled out here. We shall mention details
later when we need them.}. Thus we introduce the 
following notation to cover both situations simultaneously:
\be \label{defs1}
\begin{array}{ccccc}
  \kk & \; := \; & \QQ   & \mbox{ or } & \QQ(\rf)     \\
  \oo & \; := \; & \ZZ   & \mbox{ or } & \ZZ[\tau]    \\
  \OO & \; := \; & \JJ   & \mbox{ or } & \II          \\
  \LL & \; := \; & \ZZ^4 & \mbox{ or } & \ZZ[\tau]^4 
\end{array}
\ee
By $\LL$ we really mean the ring 
$\oo 1 + \oo \bm{i} + \oo \bm{j} + \oo \bm{k} \subset \OO$, and
we observe that
\be \label{inclusion}
   2 \OO \; \subset \; \LL \; \subset \; \OO \, .
\ee

Let us first note that $\OO$ is a maximal order in $\HH(\kk)$ 
\cite{Vigneras}. As such, each prime ideal $\PPP$ of $\OO$ corresponds
to one prime ideal of $\oo$, namely to $\ppp := \PPP\cap\oo$,
and this sets up a one-to-one correspondence between their prime
ideals \cite[Thm.\ 22.4]{Reiner}.
{}Furthermore, we have unique factorization of each 2-sided ideal
${\mathfrak A}$ of $\OO$: 
\be 
   {\mathfrak A} \; = \; \PPP_1^{k_1}\cdot\ldots\cdot\PPP_r^{k_r}
\ee
where $\PPP_1^{k_1},\ldots,\PPP_r^{k_r}$ are prime ideals and all
$k_i\geq 0$.

The prime 2 (which is prime both in $\ZZ$ and in $\ZZ[\tau]$) plays a
special role in this paper. In the case of $\II$, $2\II$ is the prime
ideal of $\II$ lying over $2\ZZ[\tau]$. The case of $\JJ$, however, is
more complicated. Here, $(1+\bm{i})\JJ=\JJ(1+\bm{i})$ is the prime
ideal lying over $2\ZZ$ and $(1+\bm{i})^2\JJ=2\JJ$ \cite{Hurwitz}.
It is this ramification of $2\ZZ$ in $\JJ$ that accounts for the
stray factor in Eq.~(\ref{Hurwitz-zeta}) and, later on, in 
Eq.~(\ref{genfun1}). To cope with this, we shall call $\bm{a}\in\JJ$ 
an {\em odd} element of $\JJ$ if $\bm{a}\not\in(1+\bm{i})\JJ$.
It is useful to note that $\bm{a}\in\JJ$ is odd if and only if 
$|\bm{a}|^2\in\ZZ$ (its quaternionic norm) is odd.

Let $\bm{a}\in\OO$. We define
\be \label{defs2}
   A(\bm{a}) \; := \; \{ r\in \kk\mid r\bm{a}\in\OO \} \, .
\ee
{}For $\bm{a}\neq 0$, $\oo\subset A(\bm{a})\subset \OO\bm{a}^{-1}
\cap \kk$. Thus $A(\bm{a})$ is a finitely generated $\oo$-module 
and hence is a fractional ideal of $\oo$ \cite[Ch.\ I.4]{Janusz}.
It follows that $A(\bm{a})^{-1}$ is an ordinary ideal of $\oo$ and,
consequently, $A(\bm{a})^{-1} = \oo c$ for some $c\in\oo$. We call
$c = c(\bm{a}) \in \oo$ (which is determined up to a unit of $\oo$) 
the {\em content} of $\bm{a}$: $A(\bm{a}) = \oo c(\bm{a})^{-1}$.
\begin{defin}
  We say that $\bm{a}\in\OO$ is {\em $\OO$-primitive} if
  $A(\bm{a})=\oo$.
\end{defin}
We have just seen that $\bm{a}$ is $\OO$-primitive if and only if
$c(\bm{a})$ is a unit in $\oo$. So, we have
\begin{lemma} \label{primitive}
 {}For any $\bm{a}\in\OO\backslash\{0\}$, $c(\bm{a})^{-1} \bm{a}$
 is a primitive element of $\OO$. \hfill $\square$
\end{lemma}

Let us now come to the link between submodules and similarities.
\begin{defin}
 Let $L$ be a $\ZZ$-module in $\RR^4$ that spans $\RR^4$.
 $\MM\subset L$ is called a {\em similarity submodule} (SSM) of $L$
 if there is an $\alpha\in\RR\backslash\{0\}$ and an $R\in{\it O}(4)$
 so that $\MM=\alpha R(\OO)\subset L$. If $L$ is a {\em lattice},
 we call $\MM$ a {\em similarity sublattice} (SSL) of $L$.
\end{defin}
Consider the case when $\MM$ is an SSM of $\OO$.
It is immediate that such an $\MM$ is an $\oo$-submodule of $\OO$.
Since $\oo$ is a principal ideal domain and $\OO$ is a free $\oo$-module
of rank 4, we see that $\MM$ is also a free $\oo$-module of rank 4.
Consequently, the index $[\OO:\MM]$ of $\MM$ in $\OO$ is finite.

We now come to the first crucial assertion in our classification of
the similarity submodules according to their indices. 
\begin{prop} \label{prop1}
 Let $\MM\subset\OO$ be a similarity submodule. Then there exist
 $\bm{a},\bm{b}\in\OO$, with $\bm{a}$ primitive, such that
 $\MM=\bm{a}\OO\bm{b}$. In addition, in the case $\OO=\JJ$, we
 can arrange for $\bm{a}$ to be odd.
\end{prop}
{\sc Proof}: By assumption, $\MM=\alpha R(\OO)$ with $\alpha\in\RR$,
$\alpha\neq 0$, and $R\in\mbox{O}(4)=\mbox{O}(4,\RR)$. Using (\ref{action}),
and noting that $\overline{\OO}=\OO$, we can write $\MM=\bm{a}\OO\bm{b}$ 
where $\bm{a}=(a^{}_0,a^{}_1,a^{}_2,a^{}_3)$, 
$\bm{b}=(b^{}_0,b^{}_1,b^{}_2,b^{}_3)$ and $\bm{a},\bm{b} \in \HH(\RR)$.

Since $2\bm{a} \{ 1,\bm{i},\bm{j},\bm{k} \} \bm{b} 
\subset 2\bm{a}\OO\bm{b}\subset\LL$, we can infer from the explicit
matrix form (\ref{matrix}) that all the matrix entries of
$2 M(\bm{a},\bm{b})$ lie in $\oo$. Combining suitable entries,
four at a time, we can obtain that
\be \label{proof1}
          8 a^{}_i b^{}_j \in \oo \, , \quad \quad
          \mbox{ for all } i,j \in \{0,1,2,3\} \, .
\ee
Since $\bm{a}\OO\bm{b}=\bm{a} r \OO r^{-1} \bm{b}$ for all 
$r\in \RR\backslash\{0\}$,
we can arrange that some $a^{}_i\in \kk\backslash\{0\}$, whereupon we
see, via (\ref{proof1}), that we can choose $\bm{a},\bm{b}\in\HH(\kk)$
without loss of generality. Clearing denominators and using 
Lemma~\ref{primitive}, we may further assume that $\MM=\bm{a}\OO\bm{b}$
with $\bm{a}\in\OO$ and $\bm{a}$ primitive.
{}From (\ref{proof1}), we now get $8b^{}_j\bm{a}\in\OO$, whence
$8b^{}_j\in A(\bm{a})=\oo$ for all $j\in\{0,1,2,3\}$. We conclude
that $8\bm{b}\in\OO$.

We now have to dispose of the factor 8 in order to prove the Proposition.
Consider the 2-sided ideal $\OO\bm{a}\OO\bm{b}\OO$ of $\OO$. Since $\OO$
is a maximal order \cite{Hurwitz,Reiner,Vigneras} in $\HH(\kk)$, we
have a unique factorization
\be \label{proof2}
    \OO\bm{a}\OO\bm{b}\OO \; = \;
    \PPP_1^{k_1}\cdot\ldots\cdot\PPP_r^{k_r}\, , 
    \quad k_1,\ldots,k_r \in \NN_0 \, ,
\ee
where $\PPP_1^{k_1},\ldots,\PPP_r^{k_r}$ are prime ideals of $\OO$
and $\PPP^{}_i\cap\oo=\ppp^{}_i$ are in one-to-one correspondence with
distinct prime ideals of $\oo$ \cite{Reiner}. Similarly, we have
$\OO\bm{a}\OO = \prod \PPP_i^{m_i}$ and 
$\OO\bm{b}\OO = \prod \PPP_i^{n_i}$, with
$k_i = m_i+n_i$ for all $i\in\{1,\ldots,r\}$, $m_i\in\NN_0$, $n_i\in\ZZ$.

Since $8\OO\bm{b}\OO\subset\OO$, the only primes for which $n_i\leq 0$
is possible are those lying over $2\in\oo$. Recall that
there is exactly one prime ideal of $\OO$ that corresponds to 2, 
and this is $\PPP^{}_1=\bm{x}\OO=\OO\bm{x}$, where $\bm{x}=(1+\bm{i})$ in the
case $\OO=\JJ$ and $\bm{x}=2$ in the case $\OO=\II$.

Thus we have $n_i\geq 0$ for all $i>1$. If we can now show that
$n_1=0$, then $\bm{b}\in\OO$ and our assertion follows.

Suppose $n_1<0$. Then $m_1\geq|n_1|>0$, and we can write
\be
   \bm{a}\OO\bm{b} \; = \; 
       \bm{a}\bm{x}^{-|n_1|}\OO\bm{x}^{|n_1|}\bm{b}
\ee
and $\OO\bm{a}\bm{x}^{-|n_1|}\OO\subset\PPP_1^{m_1-|n_1|}
\PPP_2^{m_2}\cdot\ldots\cdot\PPP_r^{m_r}\subset\OO$.
Similarly, $\OO\bm{x}^{|n_1|}\bm{b}\OO\subset\OO$.

In the icosian case (where $\bm{x}=2\in\oo$), the primitivity of
$\bm{a}$ rules out that $\bm{a}\bm{x}^{-|n_1|}\in\OO$, so $n_1<0$
is impossible here. In the Hurwitzian case, since $(1+\bm{i})^2=2$
(up to units), $n_1=-1$ is still possible. Then, 
$\bm{a}\bm{x}^{-1},\bm{xb}\in\OO$ and $\bm{a}\bm{x}^{-1}$ is still 
primitive. We take these as the new $\bm{a},\bm{b}$.
This achieves the correct form, $\MM=\bm{a}\OO\bm{b}$,
with $\bm{a},\bm{b}\in\OO$ and $\bm{a}$ $\OO$-primitive.
When $\OO=\JJ$ and $\bm{a}$ is even, we have 
$\bm{a}\in(1+\bm{i})\OO\backslash2\OO$, and we may replace
$\bm{a}\OO\bm{b}$ by $\bm{a}\bm{x}^{-1}\OO\bm{x}\bm{b}$.
\hfill $\square$

{\bf Remark 1}: Given 
$\MM=\bm{a}^{}_1\OO\bm{b}^{}_1\subset\OO$ where 
$\bm{a}^{}_1,\bm{b}^{}_1\in\HH(\kk)$, the above argument shows that
we can constructively adjust $\bm{a}^{}_1,\bm{b}^{}_1$ to
$\bm{a}=\bm{a}^{}_1 r^{-1}$, $\bm{b}=r \bm{b}^{}_1$, where $r\in\kk$
for $\OO=\II$ and $r\in\kk\cup\kk(1+\bm{i})$ for $\OO=\JJ$, so as to have
$\MM=\bm{a}\OO\bm{b}$ with $\bm{a},\bm{b}\in\OO$ and $\bm{a}$ primitive
(and odd for $\OO=\JJ$). The argument also shows that once $\bm{a}$ is
adjusted to be primitive (and odd if $\OO=\JJ$) then $\bm{b}$
necessarily lies in $\OO$.

Since we are ultimately interested in the similarity submodules, and not
so much in the actual self-similarities themselves, we have to draw the
attention now to the symmetries of our maximal orders.
The group of units $\OO^{\times}$ of $\OO$ is, geometrically, a finite
root system $\Delta$ of type\footnote{We are using the symbol $D_4$
both for the root system and for the corresponding root lattice.
The convex hull of the 24 roots of $D_4$ is the regular 24-cell
\cite[p.\ 292]{Cox73}, mentioned in an earlier footnote.
RVM would like to take this opportunity to note that in \cite{CMP}
the root system $D_4$ was inexplicably left out of the classification
of subroot systems of the root system of type $H_4$. All such $D_4$
subroot systems are conjugate by the Weyl group of $H_4$.}  
$D_4$ (resp.\ $H_4$) for $\JJ$ (resp.\ $\II$).
Any self-similarity of $\OO$ that is {\em surjective} is necessarily an
isometry (norm preserving) and so must map $\Delta$ onto itself. 
Conversely, any isometry which stabilizes $\Delta$ will also stabilize 
its $\oo$-span which is $\OO$. Thus
\be \label{stab}
    \mbox{stab}^{}_{{\rm O}(4)} \OO \; = \; \mbox{Aut}(\Delta) \, .
\ee

{}For $\OO=\II$, $\mbox{Aut}(\Delta)$ is the Weyl group $W(\Delta)$
of $\Delta$, which is $H_4$, and consequently all elements of 
$\mbox{Aut}(\Delta)^{+}$ (the orientation preserving part of
$\mbox{Aut}(\Delta)$) are realized by mappings $M(\bm{u},\bm{v})$
with $\bm{u},\bm{v} \in \II^{\times}$, i.e.\ units.

{}For $\OO=\JJ$, $[\mbox{Aut}(\Delta):W(\Delta)]=6$, the additional
symmetry being due to the diagram automorphisms of $D_4$. In fact,
$\mbox{Aut}(\Delta)$ is the Weyl group of $F_4$, and the root system
of type $F_4$ can be realized explicitly as
\be
    \tilde{\Delta} \; := \; \Delta \cup
       \{ \pm\bm{u}\pm\bm{v}\mid\bm{u},\bm{v}\in
          \{1,\bm{i},\bm{j},\bm{k}\}, \; 
          \bm{u}\neq\bm{v} \} \, ,
\ee
i.e.\ by adjoining to $\Delta$ the elements of $\JJ$ of square length 2.
This time, $\mbox{Aut}(\Delta)^{+}$ is realized as the set of mappings
$M(\bm{u},\bm{v})$ with either $\bm{u},\bm{v}\in\Delta$ or
$\bm{u},\bm{v}\in(\tilde{\Delta}\backslash\Delta)/\sqrt{2}$.
It will be observed that all the elements of $\tilde{\Delta}\backslash\Delta$
lie in the ideal $(1+\bm{i})\JJ$ (in fact they are all the generators of
this ideal). So, we have proved
\begin{prop} \label{symmetries}
 The orientation preserving self-similarities of $\OO$ onto itself are
 precisely the maps $M(\bm{u},\bm{v})$ for $\bm{u},\bm{v}\in\OO$
 with $|\bm{u}|=|\bm{v}|=1$ and, in the case $\OO=\JJ$, also the maps
 ${1\over 2}M(\bm{u},\bm{v})$ for $\bm{u},\bm{v}\in\OO$ with
 $|\bm{u}|^2=|\bm{v}|^2=2$. \hfill $\square$
\end{prop}
We say that an SSM $\bm{a}\OO\bm{b}$ is given in {\em canonical form}
if $\bm{a},\bm{b}\in\OO$ with $\bm{a}$ being $\OO$-primitive
(and with $\bm{a}$ odd if $\OO=\JJ$).
\begin{prop} \label{comp}
 Similarity submodules $\bm{a}^{}_1 \OO \bm{b}^{}_1$ and
 $\bm{a}^{}_2 \OO \bm{b}^{}_2$ written in canonical form are equal
 if and only if both
 $\bm{a}_1^{-1}\bm{a}^{}_2$ and $\bm{b}_2^{}\bm{b}^{-1}_1$
 are units in $\OO$.
\end{prop}
{\sc Proof}: Suppose that
$\bm{a}^{}_1 \OO \bm{b}^{}_1=\bm{a}^{}_2 \OO \bm{b}^{}_2$. Then
$\bm{a}^{-1}_1\bm{a}^{}_2\OO\bm{b}^{}_2\bm{b}^{-1}_1=\OO$.
According to Prop.~\ref{symmetries}, one of two things may happen:
\newline
(i) There is an $r\in\kk$ so that 
    $r\bm{a}^{-1}_1\bm{a}^{}_2=\bm{u}\in\OO^{\times}$ (and 
    $r^{-1}\bm{b}^{}_2\bm{b}^{-1}_1=\bm{v}\in\OO^{\times}$).
    Then, $r\bm{a}^{}_2=\bm{a}^{}_1\bm{u}\in\OO$ gives $r\in\oo$
    since $\bm{a}^{}_2$ is primitive. Likewise,
    $r^{-1}\bm{a}^{}_1=\bm{a}^{}_2\bm{u}^{-1}_{}\in\OO$ gives
    $r^{-1}\in\oo$, so $r\in\oo^{\times}\subset\OO^{\times}$
    and we are done in this case.
\newline
(ii) We are in the case $\OO=\JJ$ and there is an $r\in\kk$ so that
    $2r\bm{a}^{-1}_1\bm{a}^{}_2=(1+\bm{i})\bm{u}$, $\bm{u}\in\OO^{\times}$.
    This gives $2r\bm{a}^{}_2\in\OO$ and
    $r^{-1}\bm{a}^{}_1=\bm{a}^{}_2\bm{u}_{}^{-1}(1-\bm{i}) \in\OO$
    whence $2r,r^{-1}\in\ZZ$. Thus $r=\pm1,\pm{1\over 2}$.
    Since $\bm{a}^{}_1$ and $\bm{a}^{}_2$ are both odd and
    $2\in(1+\bm{i})^2 \OO^{\times}$, none of these values of $r$ is possible.
\newline
The reverse direction is clear. \hfill $\square$

We are now in the position to formulate the main result of
this section.
\begin{theorem} \label{commonthm}
 The number of similarity submodules of $\OO$ of a given index is a
 multiplicative arithmetic function. Its Dirichlet series generating
 function is given by
\be
   {}F^{}_{\OO}(s) \; = \; \frac{(\zeta^{}_{\OO}(s))^2}
                              {\zeta^{}_{\kk}(4s)} \cdot
       \cases{\frac{1}{1+4^{-s}} , & if $\OO=\JJ$ \cr
                1 ,                & if $\OO=\II$ \, . }
\ee
\end{theorem}
{\sc Proof}: As a result of Prop.~\ref{symmetries} and Prop.~\ref{comp},
any SSM $\MM$ of $\OO$ can be uniquely written as
\be \label{fac}
      \MM \; = \; \bm{a} \OO \OO \bm{b} \, ,
\ee
i.e.\ as a product of a right and a left ideal of $\OO$, where
$\bm{a}$ is $\OO$-primitive (and also odd if $\OO=\JJ$).

Any right ideal $\bm{d}\OO$ can be written {\em uniquely} as a product
$c(\bm{d}) \oo \bm{a} \OO$, where $c(\bm{d})$ (the content of $\bm{d}$) 
is in $\oo$ and $\bm{a}\in\OO$ is $\OO$-primitive. In addition, we have
the formula
$[\OO : \bm{d}\OO] = [\oo:c(\bm{d})\oo]^4 \cdot [\OO:\bm{a}\OO]$.
Thus the Dirichlet series for the {\em primitive} right ideals of $\OO$
(those $\bm{a}\OO$ with $\bm{a}$ primitive) is the quotient of
two zeta functions, 
\be \label{prim}
    \zeta^{}_{\OO}(s)/\zeta^{}_{\kk}(4s) \, .
\ee
In the case $\OO=\II$, the factorization (\ref{fac}) leads at once
to $F^{}_{\II}(s)=(\zeta^{}_{\II}(s))^2/\zeta^{}_{\QQ(\tau)}(4s)$.

In the case $\OO=\JJ$, writing (\ref{prim}) explicitly as an Euler
product, we find that the contribution for the prime 2 is
$(1+4^{-s})$. This corresponds to the fact that the primitive right
ideals are either of the form $\bm{a}\JJ$ with $\bm{a}$ odd or of the form
$\bm{a}(1+\bm{i})\JJ$ with $\bm{a}$ odd. The latter ones are to be
removed from the counting (because $(1+\bm{i})\JJ=\JJ(1+\bm{i})$ otherwise
leads to doubly counting them), and the corresponding Dirichlet series
is obtained by removing the factor $(1+4^{-s})$.
This gives the result claimed. \hfill $\square$

\subsection*{Results: lattices}

Let us now consider the Hurwitzian case $\OO=\JJ$ in detail, and also the 
lattice $\LL=\ZZ^4$. If $\bm{x}=(1+\bm{i})$, we have $\bm{x}\OO=\OO\bm{x}$
and also $\bm{x}\LL=\LL\bm{x}$. These lattices are related by
\be \label{subset1}
   \bm{x}\OO \; \stackrel{2}{\subset} \; \LL
             \; \stackrel{2}{\subset} \; \OO
\ee
where the integer on top of the inclusion symbol is the corresponding index.
\begin{lemma} \label{back}
 If $\bm{a}\LL\bm{b}\subset\LL$, then there exist $\bm{a}^{}_1,\bm{b}^{}_1
 \in\OO$, with $\bm{a}^{}_1$ odd and $\OO$-primitive,
 such that $\bm{a}\LL\bm{b}=\bm{a}^{}_1\LL\bm{b}^{}_1$. 
\end{lemma}
{\sc Proof}: Let $\bm{a}\LL\bm{b}\subset\LL$. We can assume, without loss
of generality, that $\bm{a},\bm{b}\in\HH(\QQ)$ and, by using a suitable
scaling and the fact that $\bm{x}\LL=\LL\bm{x}$, we may even assume that
$\bm{a}\in\OO$ and that $\bm{a}$ is $\OO$-primitive and odd. Since 
$\bm{a}\OO\bm{x}\bm{b}\subset\bm{a}\LL\bm{b}\subset\OO$, the conditions
on $\bm{a}$ already show that $\bm{x}\bm{b}\in\OO$ (see Remark 1). 
Write $\bm{b}=\bm{x}^{-1} \bm{c}$ with $\bm{c}\in\OO$.

Consider $[\LL:\bm{a}\LL\bm{b}]= |\bm{a}|^4 \, |\bm{b}|^4 =
\frac{|\bm{a}|^4 \, |\bm{c}|^4}{4} \in\ZZ$. Since $\bm{a}$ is odd,
it follows that $4 \mid |\bm{c}|^2$ and hence $\bm{x}\mid\bm{c}$,
because any even element in $\JJ$ is of the form $(1+\bm{i})^r$
times an odd element \cite{Hurwitz}.
Consequently, $\bm{b}\in\OO$ and we conclude that $\bm{a}\LL\bm{b}\subset\LL$
implies that we can rearrange the quaternions in the way claimed.
\hfill $\square$

This provides a link between the SSL problems for $\LL$ and for $\OO$. 
The difference between the counting arises as follows.
The symmetry group of $\OO$ is isomorphic with the Weyl group of $F_4$
(see above), while that of $\LL$, which is the Weyl group of $B_4$,
is a subgroup of index 3. As we shall show, this only influences
the number of SSLs of even index when going from $\OO$ to $\LL$.
\begin{theorem} \label{thm-lattices}
 The possible indices of similarity sublattices of hypercubic lattices
 in 4D are precisely the squares of rational integers. The number of SSLs of
 given index is a multiplicative arithmetic function. 
 {}For the case of $\JJ=D_4^*$, the corresponding Dirichlet
 series generating function $F^{}_{\JJ}$ reads
\be \label{genfun1}
   {}F^{}_{\JJ}(s) \; = \; \frac{(\zeta^{}_{\JJ}(s))^2}{(1+4^{-s})\,\zeta(4s)}
                 \; = \; \frac{(1-2^{1-2s})^2}{1+4^{-s}} \cdot
                 \frac{(\zeta(2s)\zeta(2s-1))^2}{\zeta(4s)} \, .
\ee
 The same series also applies to the lattice $D_4$,
 while for the primitive hypercubic lattice, $\ZZ^4$, it reads
\be \label{genfun2}
   {}F^{}_{\ZZ^4}(s) \; = \; (1+\frac{2}{4_{}^s}) \cdot {}F^{}_{\JJ}(s) \, .
\ee
\end{theorem}
{\sc Proof}: The statement about $F^{}_{\JJ}(s)$ follows directly from
Theorem \ref{commonthm} and from Eq.~(\ref{Hurwitz-zeta}). It is rather 
easy to see from the Euler product representation that precisely all 
squares of integers occur as indices.

In order to extend this to $\LL$, we have to understand how the different
symmetries lead to different countings. Assume $\bm{a}\LL\bm{b}\subset\LL$.
Due to Lemma~\ref{back}, we may assume that $\bm{a},\bm{b}\in\OO$ with
$\bm{a}$ odd and $\OO$-primitive, i.e.\ we assume canonical form.

By Prop.~\ref{comp}, $\bm{a}\OO\bm{b}=\bm{a}^{}_1\OO\bm{b}^{}_1$ if and only 
if there are units $\bm{u},\bm{v}\in\OO^{\times}$ with 
$\bm{a}^{}_1=\bm{a}\bm{u}$ and $\bm{b}^{}_1=\bm{v}\bm{b}$. 
However, the unit group of $\LL$ is
only the quaternion group, $Q=\{\pm 1,\pm\bm{i},\pm\bm{j},\pm\bm{k}\}$,
and $Q$ is a normal subgroup of $\OO^{\times}$ with
$\OO^{\times}/Q\simeq\ZZ/3\ZZ$.
We may take $\bm{t}=(1,1,1,1)/2$ as a suitable representative of
a generator of this cyclic group in $\OO^{\times}$ \cite[{\S} 26]{duval}.
Note that, for $\bm{u},\bm{v}\in\OO^{\times}$, we have
$\bm{a}\bm{u}\LL\bm{v}\bm{b}=\bm{a}\LL\bm{b}$ if and only if
$\bm{uv}\in Q$.
So the single SSL $\bm{a}\OO\bm{b}$ of $\OO$ {\em may} give rise to three
different SSLs of $\LL$, namely to  $\bm{a}\LL\bm{b}$,
$\bm{a}\LL\bm{t}\bm{b}$, $\bm{a}\LL\bm{t}^2\bm{b}$. Whether or not
this happens depends on whether or not the latter two are actually in $\LL$.

Now, if $\bm{a}\LL\bm{b}$ and $\bm{a}\LL\bm{t}\bm{b}$ are both in $\LL$,
then so is $\bm{a}(\LL+\LL\bm{t})\bm{b}=\bm{a}\OO\bm{b}$, whence also
$\bm{a}\LL\bm{t}^2\bm{b}\subset\LL$. Similarly,
$\bm{a}\LL\bm{t}^2\bm{b}\subset\LL$ implies 
$\bm{a}\LL\bm{t}\bm{b}\subset\LL$, too. Thus $\bm{a}\OO\bm{b}\subset\OO$
gives rise to 3 different SSLs of $\LL$ if and only if
$\bm{a}\OO\bm{b}\subset\LL$.

However, $\bm{a}\OO\bm{b}\subset\LL$ implies that
$[\OO:\bm{a}\OO\bm{b}]$ is even because $[\OO:\LL]=2$.
Conversely, $[\OO:\bm{a}\OO\bm{b}]$ even implies that
$|\bm{a}|^4\,|\bm{b}|^4$ is divisible by 4, so $\bm{a}$ or $\bm{b}$
must be even and hence $\bm{a}\OO\bm{b}\subset\bm{x}\OO\subset\LL$.
In short, the 3 SSL situation occurs if and only if
$[\OO:\bm{a}\OO\bm{b}]$ is even. 
As a consequence, the counting function for $\LL$ is still
multiplicative, and the modification in the Euler product expansion
occurs only in the factor that belongs to the prime 2. It is easy
to check that the result is that given in the Theorem.
\hfill $\square$

If we take into account that the possible indices are always squares,
it is reasonable to define the appropriate coefficients as follows,
\be \label{newseries}
     {}F^{}_{\LL}(s) \; = \; \sum_{m=1}^{\infty}
                     \frac{f^{}_{\LL}(m)}{m^{2s}} \, .
\ee
So, the coefficients actually are
\be
     f^{}_{\LL}(m) \; = \; |\{\LL' \mbox{ is SSL of } \LL
                \mid  [\LL : \LL'] = m^2 \} | \, .
\ee
To simplify explicit formulas here and later on, we introduce the function
\be \label{help}
   g(n,r) \; := \; 
    (r+1) n^r + 2 \frac{1-(r+1) n^r + r\, n^{r+1}}
                       {(n-1)^2}
\ee
for integers $r\geq 0$ and $n>1$. Note that $g(n,0)=1$.
An explicit expansion of the Euler 
factors now gives the following result.
\begin{coro} The arithmetic function $f^{}_{\JJ}(m)$ is multiplicative. 
It is given by 
\be \label{afun1}
     f^{}_{\JJ}(p^r) \; = \; 
     \cases{ 1 , & if $p=2$ \cr
            g(p,r), & if $p$ is an odd prime }
\ee
where $r\geq 0$ in all cases. Similarly, $f^{}_{\ZZ^4}(m)$ is a 
multiplicative arithmetic function. It is related to $f^{}_{\JJ}(m)$ via
\be \label{modif}
      f^{}_{\ZZ^4}(m) \; = \; 
            \cases{f^{}_{\JJ}(m) , & $m$ odd \cr
                   3 \cdot f^{}_{\JJ}(m) , & $m$ even.} 
\ee
\end{coro}
The first few terms of $F^{}_{\JJ}(s)$ read explicitly
\be
  {}F^{}_{\JJ}(s) \; = \; 
 \mbox{\small $
     1+\frac{1}{4^s}+\frac{8}{9^s}+\frac{1}{16^s}+
       \frac{12}{25^s}+\frac{8}{36^s}+\frac{16}{49^s}+
       \frac{1}{64^s}+\frac{41}{81^s}+\frac{12}{100^s}+
       \frac{24}{121^s}+\frac{8}{144^s}+  \cdots $ }
\ee
while those for $F^{}_{\ZZ^4}(s)$ follow easily from (\ref{modif}). They are now 
listed as \cite[sequence A 045771]{Sloane} and \cite[sequence A 035292]{Sloane}, 
respectively. Note that Eq.~(\ref{afun1}) implies that the SSMs of $\JJ$ of index 
$4^r$ are unique -- they are, in fact, just the 2-sided ideals $(1+\bm{i})^r \JJ$.

Let us now, in line with the previous examples, briefly consider the
asymptotic behaviour of the coefficients. Since $\zeta(s)\neq 0$
in $\{\mbox{Re}(s)\geq 1\}$, it is clear that $F^{}_{\LL}(s)$ is
meromorphic in the same half-plane, with only one pole which is 
of second order and located at $s=1$. 
This is true both of $\LL=\JJ$ and $\LL=\ZZ^4$.
Using the Dirichlet series (\ref{newseries}) and applying the results 
from the Appendix to $F^{}_{\LL}(s/2)$, we get the following
\begin{coro}
The coefficients $f^{}_{\LL}(m)$ grow faster than linear on average for 
large $m$, and we have the asymptotic behaviour
\be
   \sum_{m\leq x} f^{}_{\LL}(m) \; \simeq \; C^{}_{\LL} \cdot x^2 \log(x)
   \quad \quad \quad \mbox{(as $x\to\infty$)}
\ee
where the constant is given by
\be
   C^{}_{\LL} \; = \; \frac{1}{2} \, 
              \mbox{\rm res}_{s=1} \, ((s-1) {}F^{}_{\LL}(s))
   \; = \; \frac{1}{4} \cdot 
           \cases{ 1   \, , & if $\LL=\JJ$ \cr
           \frac{3}{2} \, , & if $\LL=\ZZ^4$ \, . } 
\ee
\end{coro}
Note that this really is an asymptotic result, and that, numerically,
the estimates of $C^{}_{\LL}$ converge rather slowly (from above) to the
values given in the Corollary.

\subsection*{Results: modules}

Let us first state the result for the icosian ring $\II$ itself.
\begin{theorem} \label{thm-modules}
 The possible indices of similarity submodules of the icosian ring
 $\II$ are the squares of rational integers that can be represented
 by the quadratic form $x^2 + xy - y^2$. The number of SSMs of
 given index is a multiplicative arithmetic function, its 
 Dirichlet series generating function reads
\be \label{genfun3}
   {}F^{}_{\II}(s) \; = \; \frac{(\zeta^{}_{\II}(s))^2}
                              {\zeta^{}_{\kk}(4s)}
                 \; = \; 
                 \frac{(\zeta^{}_{\kk}(2s)\zeta^{}_{\kk}(2s-1))^2}
                      {\zeta^{}_{\kk}(4s)}
\ee
where $\kk=\QQ(\tau)$.
\end{theorem}
The proof follows immediately from Theorem \ref{commonthm} in combination
with Eq.~(\ref{icosian-zeta}). Note that the possible indices are just
the squares of the possible norms of ideals in $\QQ(\tau)$ and hence
of the form given. Taking this into account, we write
\be
   {}F^{}_{\II}(s) \; = \; \sum_{m=1}^{\infty}
                 \frac{f^{}_{\II}(m)}{m^{2s}}
\ee
in analogy to above, and obtain, by an explicit expansion of the Euler 
factors, the following result (compare \cite[sequence A 035284]{Sloane}).
\begin{coro} The arithmetic function $f^{}_{\II}(m)$ is multiplicative
 and given by
\be
   f^{}_{\II}(p^r) \; = \; 
    \cases{ g(5,r) ,     & if $p=5$ \cr
            0 ,          & if $p\equiv\pm 2$ mod 5 and $r$ is odd \cr
            g(p^2,\ell), & if $p\equiv\pm 2$ mod 5 and $r=2\ell$ \cr
            \sum_{s=0}^{r} g(p,s) g(p,r\!-\!s) , &
                        if $p\equiv\pm 1$ mod 5 }
\ee
where always $r\geq 0$ and $g$ is the function defined in Eq.~(\ref{help}).
\end{coro}
The first few terms of $F^{}_{\II}(s)$ read explicitly
\be
  {}F^{}_{\II}(s) \; = \; 
 \mbox{\small $
     1+\frac{10}{16^s}+\frac{12}{25^s}+\frac{20}{81^s}+
       \frac{48}{121^s}+\frac{66}{256^s}+\frac{80}{361^s}+
       \frac{120}{400^s}+\frac{97}{625^s}+\frac{120}{841^s}+
       \frac{128}{961^s}+  \cdots $ }
\ee
Let us briefly look at the $10$ SSMs of $\II$ of index $16$. {}From 
Eq.~(\ref{ico-zeta-2}) it is obvious that they are just the 5 left ideals 
$\II\bm{a}$ and the 5 right ideals $\bm{a}\II$, with suitable generators 
$\bm{a}$ with $\mbox{N}(|\bm{a}|^2)=16$. Note that none of them is 2-sided.

{}Finally, we can again determine the asymptotic behaviour along the
lines used before. $F^{}_{\II}(s/2)$ is holomorphic in the half-plane
$\{\mbox{Re}(s)\geq 2 \}$, with a single second-order pole at $s=2$. 
With the results from the Appendix, we then obtain
\begin{coro}
{}For $x\to\infty$, the asymptotic behaviour of the coefficients
$f^{}_{\II}(m)$ is
\be
   \sum_{m\leq x} f^{}_{\II}(m) 
         \; \sim \; \frac{6 (\log(\tau))^2}{5\sqrt{5}} \, x^2 \log(x)
         \; \simeq \; 0.124271 \, x^2 \log(x) \, .
\ee
\end{coro}

At this point, it would be interesting to relate these findings to the
corresponding ones for the $\ZZ[\tau]$-modules $\LL=\ZZ[\tau]^4$ and
$\MM=\JJ[\tau]=\II\cap\II'$. Since this requires a lot more effort than
in the previous case ($\ZZ^4$ versus $\JJ$), we postpone it, and rather
state the result for the cubian\footnote{The term ``octonian'' would be 
more natural a choice, but it has already been taken!} maximal order $\KK$ in 
$\HH(\QQ(\rz))$.

\begin{theorem} \label{thm-modules-2}
 The possible indices of similarity submodules of the cubian ring
 $\KK$ are the squares of rational integers that can be represented
 by the quadratic form $x^2 - 2 y^2$. The number of SSMs of
 given index is a multiplicative arithmetic function, its 
 Dirichlet series generating function reads
\begin{eqnarray} \label{genfun4}
   {}F^{}_{\KK}(s) & = & \frac{(\zeta^{}_{\KK}(s))^2}
                              {\zeta^{}_{\kk}(4s)}
                 \; = \; 
                 \frac{(\zeta^{}_{\kk}(2s)\zeta^{}_{\kk}(2s-1))^2}
                      {\zeta^{}_{\kk}(4s)}   \\  &  = &
\mbox{\small $
     1+\frac{6}{4^s}+\frac{22}{16^s}+\frac{32}{49^s}+
       \frac{66}{64^s}+\frac{20}{81^s}+\frac{192}{196^s}+
       \frac{178}{256^s}+\frac{72}{289^s}+\frac{120}{324^s}+
       \frac{96}{529^s}+{52\over625^s}+  \cdots $ }   \nonumber
\end{eqnarray}
where $\kk=\QQ(\rz)$.
\end{theorem}
The proof follows directly from the proof of Theorem \ref{commonthm}, since
literally every step taken for the icosian ring translates into one here,
with $\kk=\QQ(\rz)$, $\oo=\ZZ[\rz]$, $\OO=\KK$, and $\LL=\ZZ[\rz]^4$.
Note that now, since 2 is not a prime in $\ZZ[\rz]$ (it actually ramifies
there), $\PPP^{}_1=\rz\,\KK=(2+\rz)\KK$ is the prime ideal of $\KK$ sitting 
on top of the prime ideal $(2+\rz)\ZZ[\rz]$, and the arguments in the proofs
have to be adjusted accordingly.

One can again work out the coefficients $f^{}_{\KK}(m)$ explicitly
(see \cite[sequence A 035285]{Sloane})
\be
   f^{}_{\KK}(p^r) \; = \; 
    \cases{ g(2,r) ,     & if $p=2$ \cr
            0 ,          & if $p\equiv\pm 3$ (8) and $r$ is odd \cr
            g(p^2,\ell), & if $p\equiv\pm 3$ (8) and $r=2\ell$ \cr
            \sum_{s=0}^{r} g(p,s) g(p,r\!-\!s) , &
                        if $p\equiv\pm 1$ (8) }
\ee
and the asymptotic behaviour, using the Appendix, is
\be
  \sum_{m\leq x} f^{}_{\KK}(m) 
         \; \sim \; \frac{15 (\log(1+\rz))^2}{22\sqrt{2}} \, x^2 \log(x)
         \; \simeq \; 0.374519 \, x^2 \log(x) \, .
\ee

\subsection*{Concluding remarks}

As we have demonstrated above, the similarity submodules of certain 4D 
$\ZZ$-modules can be classified by means of algebraic methods based on 
quaternionic algebras and their maximal orders. Together with the results
of \cite{BM1,BM2}, this essentially covers the cases related to root systems
in dimensions $d\leq 4$.

Although we did not emphasize it, one can also determine the actual
semigroups of self-similarities of these modules explicitly, notably
through the canonical representation of SSMs (Prop.\ \ref{symmetries}) and
their uniqueness up to symmetries (Prop.\ \ref{comp}). We have described
this in more detail for other cases \cite{BM1}, and the interested reader
will find no difficulty to extend that approach to this situation.

One application is concerned with the symmetries of coloured versions
of the lattices and modules under consideration. Assume that $L$ has
a non-trivial (irreducible) point symmetry, and a sublattice $L'$ which is 
the image of a self-similarity of $L$ of index $m=[L:L']>1$. If we assign 
$m$ different colours to the cosets of $L'$, certain subgroups of
the point group of $L'$ (which is conjugate to that of $L$) will 
give rise to a {\em colour symmetry} in the sense that their elements induce 
a unique, global permutation of the colours, compare \cite{Schwarz2,Ron} and
references therein.


This is also closely related to the classification of coincidence site
submodules, i.e.\ of submodules that can be written as the intersection
of the original module with a rotated copy of itself, see \cite{Baake}
for background and some recent results. Here are several open questions,
particularly in spaces of even dimension, which the above results should
help to solve for dimension four.

{}Finally, one would like to know to what extent a generalization of our
results is possible. The root lattices seem to form a sufficiently well-behaved
class of objects to try, and some partial answers on the existence of 
similarity sublattices and their possible indices are given in \cite{CRS}.
We are, however, not aware of general results along the lines discussed here, 
i.e.\ including the determination of the number of SSLs of a given index,
nor even of a method to overcome the dependence on special features such as 
the arithmetic of quaternions.

\vspace{5mm}
\subsection*{Acknowledgements}

We are grateful to Alfred Weiss for his help in understanding the
arithmetic of quaternionic maximal orders.
It is our pleasure to thank Peter Pleasants and Johannes Roth 
for several helpful discussions and Neil Sloane for communication of
material prior to publication. 
This work was supported by the German Science {}Foundation (DFG) and by 
the Natural Sciences and Engineering Research Council of Canada (NSERC).

\clearpage
\subsection*{Appendix}

In what follows, we briefly summarize the results from analytic number
theory that we need to determine certain asymptotic properties of the
coefficients of Dirichlet series generating functions. {}For the general
background, we refer to \cite{Apostol} and \cite{Zagier}.

Consider a Dirichlet series of the form 
$F(s)=\sum_{m=1}^{\infty} a(m) m^{-s}$. We are mainly interested 
in the quantity $A(x)=\sum_{m\leq x} a(m)$ and its behaviour for 
large $x$. Let us give one classical result (based upon Tauberian 
theorems) for the case that $a(m)$ is real and non-negative.
\begin{theorem} \label{meanvaluetheorem}
 Let $F(s)$ be a Dirichlet series with non-negative coefficients
 which converges for ${\rm Re}(s) > \alpha > 0$. Suppose that $F(s)$
 is holomorphic at all points of the line $\{ {\rm Re}(s) = \alpha \}$
 except at $s=\alpha$. Here, when approaching $\alpha$ from
 the half-plane right of it, we assume $F(s)$ to have a singularity
 of the form $F(s) = g(s) + h(s)/(s-\alpha)^{n+1}$ where
 $n$ is a non-negative integer, and both $g(s)$ and $h(s)$ are 
 holomorphic at $s=\alpha$. Then we have, as $x\rightarrow\infty$,
\be \label{meanvalues}
     A(x) \; := \; \sum_{m\leq x} a(m)
          \; \sim \;  \frac{h(\alpha)}{\alpha\cdot n!}
             \; x_{}^{\alpha} \, (\log(x))_{}^n \, .
\ee
\end{theorem}
The proof follows easily from Delange's theorem, e.g.\ by taking
$q=0$ and $\omega=n$ in Tenenbaum's formulation of it, 
see \cite[ch.\ II.7, Thm.\ 15]{Tenenbaum} and references given there. 

Note that Delange's theorem is a lot more general in that is still
gives results for other local behaviour of $F(s)$ in the neighbourhood
of $s=\alpha$, in particular for $n$ not an integer and even for 
combinations with logarithmic singularities.
Let us also point out that there are various extensions to
Dirichlet series with complex coefficients, e.g.\ Thm.~1 on p.\ 311
of \cite{Lang}, and even stronger results (with good error estimates)
for multiplicative arithmetic functions $a(m)$ with values in the
unit disc, see \cite[ch.\ I, {\S} 3.8 and ch.\ III, {\S} 4.3]{Tenenbaum}
for details.

The critical assumption in Theorem~\ref{meanvaluetheorem} is the 
behaviour of $F(s)$ along the entire line 
$\{ \mbox{Re}(s) = \alpha \}$. In all cases that appear
in this article, this can be checked explicitly. To do so, we have
to know a few properties of the Riemann zeta function, $\zeta(s)$,
and of the Dedekind zeta functions of $\QQ(\tau)$ and $\QQ(\rz)$.
It is well known that $\zeta(s)$ is a meromorphic function in the
complex plane, and that is has a sole simple pole at $s=1$
with residue 1 \cite[Thm.\ 12.5(a)]{Apostol}. 
It has no zeros in the half-plane 
$\{ \mbox{Re}(s)\geq 1 \}$ \cite[ch.\ II.3, Thm.\ 9]{Tenenbaum}.
The values of $\zeta(s)$ at positive even integers are known 
\cite[Thm.\ 12.17]{Apostol} and we have
\be
     \zeta(2) \; = \; \frac{\pi^2}{6} \quad , \quad
     \zeta(4) \; = \; \frac{\pi^4}{90} \, .
\ee
This is all we need to know for this case.

The Dedekind zeta function of $\kk=\QQ(\tau)$ has some similarly nice
properties. It follows from \cite[Thm.\ 4.3]{Wash} or from 
\cite[\S 11, Eq.\ (10)]{Zagier} that it can be written as
\be
     \zeta^{}_{\QQ(\tau)}(s) \; = \; \zeta(s) \cdot L(s,\chi)
\ee
where $L(s,\chi)$ is the $L$-series of the primitive Dirichlet character
\cite[Ch.\ 6.8]{Apostol} $\chi$ defined by
\be
   \chi(n) \; = \; \cases{0, & $n\equiv 0$ (5) \cr
                          1, & $n\equiv\pm 1$ (5) \cr
                         -1, & $n\equiv\pm 2$ (5) \, . } 
\ee
Since $\chi$ is not the principal character,
$L(s,\chi)=\sum_{m=1}^{\infty}\chi(m)\, m^{-s}$ is an entire function 
\cite[Thm.\ 12.5]{Apostol}.
Consequently, $\zeta^{}_{\QQ(\tau)}(s)$ is meromorphic, and
its only pole is simple and located at $s=1$. The residue is
$L(1,\chi)$ and from \cite[Thm.\ 4.9]{Wash} we get
\be
    \mbox{res}_{s=1} \, \zeta^{}_{\QQ(\tau)}(s)
    \; = \; L(1,\chi)
    \; = \; \frac{2\log(\tau)}{\sqrt{5}}
    \; \simeq \; 0.430409 \, .
\ee

Since $L(s,\chi)\neq 0$ for $\mbox{Re}(s)>1$, see \cite[p.\ 31]{Wash},
$\zeta^{}_{\QQ(\tau)}(s)$ cannot vanish there either.
Also, one can again calculate the values of $\zeta^{}_{\QQ(\tau)}(s)$
at positive even integers. This is done by means of the functional 
equation of $L(s,\chi)$
\cite[p.\ 30]{Wash} and the knowledge of the values of $L$-functions
at negative integers in terms of generalized Bernoulli numbers
\cite[Thm.\ 4.2]{Wash}, see \cite[Prop.\ 4.1]{Wash} for a formula
for them. Working this out explictly for $s=2$ and $s=4$ gives
\be
   \zeta^{}_{\QQ(\tau)}(2) 
       \; = \; \frac{2\pi^4}{75 \sqrt{5}}  
       \quad , \quad
   \zeta^{}_{\QQ(\tau)}(4)
       \; = \; \frac{4\pi^8}{16875 \sqrt{5}} \, .
\ee
Let us add that these results, and those to follow, can also be found, in rather 
explicit form, in \S 9 and \S 11 of \cite{Zagier}.

In the same way, one can determine the zeta function of $\QQ(\rz)$
and its properties. One has
$\zeta^{}_{\QQ({\scriptscriptstyle \sqrt{2}})}(s) = \zeta(s) L(s,\chi)$,
now with the primitive Dirichlet character
$\chi(n)=0,1,-1$ for $n$ even, $n\equiv\pm 1$ (8), $n\equiv\pm 3$ (8),
respectively. This zeta function has again only one simple pole, at $s=1$,
with residue $L(1,\chi)=\log(1+\rz)/\rz\simeq 0.623225$. {}Finally, we have
\be
   \zeta^{}_{\QQ({\scriptscriptstyle \sqrt{2}})}(2) 
       \; = \; \frac{\pi^4}{48 \sqrt{2}}  
       \quad , \quad
   \zeta^{}_{\QQ({\scriptscriptstyle \sqrt{2}})}(4)
       \; = \; \frac{11 \pi^8}{69120 \sqrt{2}} \, .
\ee

\end{document}